\documentclass[12pt]{amsart}
\usepackage[dvips]{graphicx}
\usepackage[usenames,dvipsnames]{color}
\usepackage{amsmath,amsthm,amsopn,amssymb,a4wide,array,enumerate,times}
\usepackage[mathscr]{eucal}
\newtheorem{theorem}{Theorem}[section]
\newtheorem{lemma}[theorem]{Lemma}

\newtheorem{proposition}[theorem]{Proposition}
\newtheorem{corollary}[theorem]{Corollary}

\renewcommand{\t}[1]{{\tilde{#1}}}

\newcommand{\ran}{{{\mathrm{ran}\,}}}
\newcommand{\clran}{{\overline{\mathrm{ran}}\,}}
\newcommand{\Ls}[1]{{{\mathbf L}({\mathcal #1}) }}

\def\card{{\mathrm{card}\,}}
\def\t#1{{\tilde{#1}}}

\def\ran{{\text{ran\,}}}

\def\clran{{\overline{\text{ran}}\,}}
\def\L#1#2{{{\bf L}({\mathcal #1,#2}) }}
\def\Ls#1{{{\mathbf L}({\mathcal #1}) }}

\def\sp#1{{{\mathcal #1} }}

\def\HardyD#1{{H^2_{\sp #1}({\mathbb D})}}

\def\HinfDL#1#2{{H^\infty_{\mathbf L(\sp #1,\sp #2)}({\mathbb D})}}
\begin{document}
\title{Outer Factorizations in One and Several Variables}
\author[M.A.~Dritschel and H.J.~Woerdeman]{
Michael A.~Dritschel$^1$
and Hugo J.~Woerdeman$^2$}
\address{School of Mathematics and Statistics\\
  Merz Court,\\ University of Newcastle upon Tyne\\
  Newcastle upon Tyne\\
  NE1 7RU\\
  UK}
\email{m.a.dritschel@newcastle.ac.uk}
\address{Department of Mathematics\\
  The College of William \& Mary\\
  Williamsburg, Virginia 23185-8795,{\rm and }\\
  Department of Mathematics\\
  Katholieke Universiteit Leuven\\
  Celestijnenlaan 200B\\
  B3001 Heverlee\\
  Belgium}
\email{hugo@math.wm.edu}
\thanks{${}^1$Research supported by the Engineering and Physical
  Sciences Research Council (EPSRC).\\
  \indent
  ${}^2$Research supported in part by National Science Foundation
  (NSF), as well as a Faculty Research Assignment (FRA) Grant from the
  College of William \& Mary.}
\abstract A multivariate version of Rosenblum's Fej\'er-Riesz theorem
on outer factorization of trigonometric polynomials with operator
coefficients is considered.  Due to a simplification of the proof of
the single variable case, new necessary and sufficient conditions for
the multivariable outer factorization problem are formulated and
proved.  \endabstract

\maketitle

\section{Introduction}
\setcounter{equation}{0} The Fej\'er-Riesz theorem for trigonometric
polynomials $q(z)=\sum_{i=-n}^nq_iz^i$ states that $q(z)\ge 0$, $z \in
{\mathbb T}$, if and only if there exists an analytic polynomial
$p(z)=\sum_{i=0}^n p_i z^i$ so that $q(z)=|p(z)|^2$, $z \in {\mathbb
  T}$.  In addition, one may choose $p$ to be void of roots inside the
open unit circle ${\mathbb D}$ (that is, $p$ is {\it outer}).  Though
simple to state and prove (use the fundamental theorem of
algebra---see, for example, \cite{MR91g:00002}),
the lemma has many useful applications; for example, in filter design,
$H^\infty$ control, and wavelet theory.  The first generalizations of
the lemma involved matrix valued trigonometric polynomials
(\cite{MR19:1098c}, \cite{MR30:1409}) and subsequently operator valued
trigonometric polynomials (in \cite{MR30:5182} a compactness condition
appears, in \cite{MR37:3378} the general operator case is done).

The present paper grew out of an interest in a multivariate analog of
the Fej\'er-Riesz theorem.  As is well-known, extensions of such
results to several variables are far from straightforward.  One of the
earliest efforts in this direction is Hilbert's well-known observation
that not all nonnegative polynomials in several real variables are
necessarily sums of squares of polynomials.  While Hilbert's result
concerns polynomials on ${\mathbb R}^d$, a similar phenomenon occurs
in the setting of trigonometric polynomials on the $d$-torus ${\mathbb
  T}^d$, where ${\mathbb T}=\{ z \in {\mathbb C} \ : \ |z|=1\}$.
Indeed, it follows from Hilbert's result (see \cite{CP} and
\cite{MR27:1779}) that a trigonometric polynomial $q(z)$ in $d$
variables of degree $(n_1,\ldots , n_d )$ that takes on nonnegative
values on ${\mathbb T}^d$, is not necessarily of the form
\begin{equation} 
  \label{sos}
  q(z) = \sum_{i=1}^k | p_i(z)|^2 , z \in {\mathbb T}^d ,
\end{equation} 
where $p_i$ are polynomials of degree $(n_1,\ldots , n_d )$.  It turns
out (see \cite{MR96f:11058}; see also \cite {parrilo},
\cite{MR2002k:90074}) that checking whether $q$ can be factored in
this way is a semidefinite feasibility problem.  In this paper, we
investigate which multivariable trigonometric polynomials are single
squares; that is, we would like $k=1$ in the above representation
\eqref{sos} with similar restrictions on the degree of the polynomial
in the factorization.

As might be expected, not putting any restrictions on the degrees of
the polynomials $p_i$ in \eqref{sos} enables more nonnegative
trigonometric polynomials to be factored.  In fact, in \cite{MAD} it
was shown that any strictly positive trigonometric polynomial (i.e.,
$q(z) >0$, $z \in{\mathbb T}^d$) allows a representation \eqref{sos}
where $p_i$ are polynomials of potentially very high degree.  This in
turn relates to factorization of real polynomials as sums of squares
of rational functions with fixed denominators \cite{MAD}.  An
important tool in \cite{MAD} is the use of Schur complements.
Inspired by this we also use the Schur complement as our main tool.
This allows for a very simple proof of Rosenblum's operator valued
Fej\'er-Riesz theorem.  The main observation in this proof is that the
sequence of finitely supported Schur complements of a banded positive
semidefinite Toeplitz operator, have a very simple inheritance
structure (see Proposition \ref{scs}).  In fact, beyond a certain
matrix size (determined by the number of nonzero diagonals) as the
Schur complement is increased one dimension in size, it is constructed
by bordering the previous Schur complement with the coefficients of
the underlying trigonometric polynomial.  Recognizing this, the task
becomes to determine the multivariate analog of this inheritance
structure.  Clearly, there are now many canonical shifts; how does one
use these?  As we will see, for the multivariate trigonometric
polynomial to have an outer factorization of the required type, the
Schur complement of the corresponding Toeplitz operator needs to
decompose in a certain way.  Subsequently, to obtain the next Schur
complement, the different terms in this decomposition need to be
shifted in different ways.  Bordering the result with the coefficients
of the trigonometric polynomial then yields the next Schur complement.
The precise statement is given in Theorem~\ref{2var}.

The paper is organized as follows.  In Section 2 we derive several
useful new properties of Schur complements.  In Section 3 we use these
newly observed properties to provide easy proofs for Rosenblum's
version of the operator valued Fej\'er-Riesz theorem and the existence
of inner-outer factorizations.  In Section 4 the multivariate case is
addressed.

\section{Auxiliary results on Schur complements}
\label{sec:lemmas}
\setcounter{equation}{0}

We will number rows and columns of an $n\times n$ matrix with
$0,\ldots,n-1$.  For $\Lambda \subseteq \{0,\ldots,n-1\}$ and an
$n\times n$ operator matrix $M$, we write $S(M;\Lambda )$ (or
$S(\Lambda)$ when there is no chance of confusion) for the Schur
complement supported on rows and columns labelled by elements of
$\Lambda$.  It is usual to view $S(\Lambda)$ as an $m\times m$ matrix,
where $m = \card \Lambda$, however it is often useful to take
$S(\Lambda)$ to be an $n\times n$ matrix.  If $\Lambda =
\{n_0,\ldots,n_{m-1}\}$, then this is done by putting the $(j,k)$
entry of $S(\Lambda)$ as an $m\times m$ matrix into the $(n_j,n_k)$
place and padding with zeros.  We use the same notation for both
versions of the Schur complement, since it should be clear from the
context which we are using.  Finally, as a further bit of notational
convenience, we write $S(M;k)$ (or $S(k)$) when $\Lambda =
\{0,\ldots,k\}$.

\begin{lemma}\label{lemma1}
  Let
  \begin{equation}
    \label{eq:1}
    M =
    \begin{pmatrix}
      A & B \\ B^* & C
    \end{pmatrix}
    =
    \begin{pmatrix}
      P^* & Q^* \\ 0 & R^*
    \end{pmatrix}
    \begin{pmatrix}
      P & 0 \\ Q & R
    \end{pmatrix} .
  \end{equation}
  Then $S(0)$ equals $P^*P$ if and only if $\ran Q \subseteq \clran
  R$.  Furthermore, for any $P$ such that $P^*P = S(0)$ and any $R$
  such that $R^*R = C$, there is a $Q$ such that (\ref{eq:1}) holds.
\end{lemma}

\begin{proof}
  Since $C = R^*R$, there is an isometry $V: \clran R \to \clran
  C^{1/2}$ such that $C^{1/2} = VR$.  Clearly $\begin{pmatrix} A & B
    \\ B^* & C \end{pmatrix} \geq 0$, so there is a contraction $G:
  \clran C^{1/2} \to \clran A^{1/2}$ with $B = A^{1/2} G C^{1/2}$, and
  consequently $B = A^{1/2} G VR$.  We also have $B = Q^*R$, so the
  assumption that $\ran Q \subseteq \clran R$ implies that $A^{1/2} G
  V = Q^*$.  Moreover, since $VV^* = 1_{\clran C^{1/2}}$, we get that
  $A^{1/2} G = Q^*V^*$.  We calculate the Schur complement
  \begin{equation*}
    \begin{split}
      S(0) &= A^{1/2} (1-GG^*) A^{1/2} = A - Q^*V^*VQ = A -
      Q^*1_{\clran R^{1/2}} Q\\ &= A-Q^*Q = P^*P+Q^*Q-Q^*Q = P^*P.\\
    \end{split}
  \end{equation*}

  Conversely assume $P^*P = S(0)$.  Then
  \begin{equation}
    \label{eq:2}
    \begin{pmatrix}
      A & B \\ B^* & C
    \end{pmatrix}
     -
     \begin{pmatrix}
       P^*P & 0 \\ 0 & 0
     \end{pmatrix}
     =
     \begin{pmatrix}
       Q^*Q & Q^*R \\ R^*Q & R^*R
     \end{pmatrix}
  \end{equation}
  If we set $V_R$, $V_Q$ to be the inclusions of $\clran R$ and
  $\clran Q$ into $\sp H$, respectively, then $Q^*R = Q^*GR$, where $G
  = V_Q^*V_R$.  By construction the Schur complement of the right side
  of (\ref{eq:2}) is zero, which implies that
  \begin{displaymath}
    0 = 1_{\clran Q}-GG^* = 1_{\clran Q}-V_Q^*V_RV_R^*V_Q =
    V_Q^*(1-P_{\clran R})V_Q,
  \end{displaymath}
  where $P_{\clran R}$ is the orthogonal projection onto $\clran R$.
  Thus $P_{\clran R} | \clran Q = 1_{\clran Q}$, and hence $\ran Q
  \subseteq \ran P_{\clran R} = \clran R$.
  
  Finally, suppose $P^*P = S(0)$ and $R^*R = C$.  Then $A - P^*P \geq
  0$ and $B = A^{1/2} G R$ for some contraction $G: \clran R \to
  \clran A^{1/2}$.  Hence $P^*P = A^{1/2}(1 - GG^*) A^{1/2}$ and there
  exists $D_G$ such that $D_GD_G^* = 1-GG^*$ and $P = A^{1/2} D_G$.
  Then setting $Q^* = A^{1/2} G$, we have $M =
  \begin{pmatrix} P^* & Q^* \\ 0 & R^* \end{pmatrix}
  \begin{pmatrix} P & 0 \\ Q & R \end{pmatrix}$, and $\ran Q \subseteq
  \ran G^* \subseteq \clran R$.
\end{proof}

\begin{lemma}\label{lemma2}
  Suppose
  \begin{equation}
    \label{eq:3}
    M =
    \begin{pmatrix}
      A&B&C \\ B^*&D&E \\ C^*&E^*&F
    \end{pmatrix}
     =
    \begin{pmatrix}
      P^*&Q^*&R^* \\ 0&S^*&T^* \\ 0&0&U^*
    \end{pmatrix}
    \begin{pmatrix}
      P&0&0 \\ Q&S&0 \\ R&T&U
    \end{pmatrix},
  \end{equation}
where $M$ is acting on ${\mathcal H}_1 \oplus {\mathcal H}_2
\oplus {\mathcal H}_3$.
  Then
  \begin{equation}
    \label{eq:4}
    S(1) - S(0)=
    \begin{pmatrix}
      Q^* \\ S^*
    \end{pmatrix}
    \begin{pmatrix}
      Q & S
    \end{pmatrix}
  \end{equation}
  if and only if
  \begin{equation}
    \label{eq:5}
    \qquad \ran Q \subseteq \clran S \qquad\text{\rm and } \qquad \ran T
    \subseteq \clran U.
  \end{equation}
  Furthermore there exists a factorization of $M$ as in (\ref{eq:3})
  with the factors operators on ${\mathcal H}_1 \oplus {\mathcal H}_2
  \oplus {\mathcal H}_3$ such that (\ref{eq:4}) and (\ref{eq:5}) hold,
  $\begin{pmatrix} P^* & Q^* \\ 0 & R^* \end{pmatrix}
  \begin{pmatrix} P & 0 \\ Q & R \end{pmatrix} = S(1)$ and $P^*P =
  S(0) = S(S(1);0)$.
\end{lemma}

\begin{proof}
  To begin with, suppose (\ref{eq:4}) holds.  Then if $ {\t P}^* \t P
  = S(0)$, we have
  \begin{equation*}
    S(1) =
    \begin{pmatrix}
      {\t P}^* & Q^* \\ 0 & S^*
    \end{pmatrix}
    \begin{pmatrix}
      {\t P} & 0 \\ Q & S
    \end{pmatrix}.
  \end{equation*}
  As $U^*U=F$, by Lemma \ref{lemma1} there exist $\t R$ and $\t T$
  such that
  \begin{equation*}
    M =
    \begin{pmatrix}
      {\t P}^*&Q^*&{\t R}^* \\ 0&S^*&{\t T}^* \\ 0&0&U^*
    \end{pmatrix}
    \begin{pmatrix}
      {\t P}&0&0 \\ Q&S&0 \\ {\t R}&{\t T}&U
    \end{pmatrix}.
  \end{equation*}
  By Lemma \ref{lemma1}
  \begin{equation*}
    \ran
    \begin{pmatrix}
      Q \\ {\t R}
    \end{pmatrix}
    \subseteq\clran
    \begin{pmatrix}
      S & 0 \\ {\t T} & U
    \end{pmatrix}
  \end{equation*}
  and
  \begin{equation*}
    \ran
    \begin{pmatrix}
      {\t R} & {\t T}
    \end{pmatrix}
    \subseteq \clran U.
  \end{equation*}
  Hence $\ran {\t T}\subseteq \clran U$ and so
  \begin{equation*}
    \clran
    \begin{pmatrix}
      S & 0 \\ {\t T} & U
    \end{pmatrix}
     = \clran S \oplus \clran U.
  \end{equation*}
  Thus $\ran Q \subseteq \clran S$.
  
  Next observe that $D = S^*S+T^*T = S^*S+{\t T}^*\t T$ and so there
  is an isometry $V_T:\clran T \to \clran {\t T}$ such that $T^* = {\t
    T}^* V_T^*$.  Also $\ran {\t T} \subseteq \clran U$ implies that
  $\ran V_T \subseteq \clran U$.  Thus $V_T$ is an isometry from
  $\clran T$ into $\clran U$.  But $U^*T = E^* = U^*{\t T} = U^*V_T
  T$, so $V_T = 1_{\clran T}$ and $\ran T \subseteq \clran U$.
  
  Now conversely assume we have a factorization of $M$ as in
  (\ref{eq:3}) where (\ref{eq:5}) holds.  Set
  \begin{equation*}
    L =
    \begin{pmatrix}
      D & E \\ E^* & F
    \end{pmatrix}
    =
    \begin{pmatrix}
      S^* & T^* \\ 0 & U^*
    \end{pmatrix}
    \begin{pmatrix}
      S & 0 \\ T & U
    \end{pmatrix}.
  \end{equation*}
  Using Lemma \ref{lemma1}, suppose $\t G$ is any other operator
  matrix satisfying ${\t G}^* {\t G}= M$ with
  \begin{equation}
    \label{eq:6}
    \t G =
    \begin{pmatrix}
      \t P&0&0 \\ \t Q&\t S& 0 \\ \t R & \t T &U
    \end{pmatrix}
  \end{equation}
  where
  \begin{equation*}
    S(1) =
    \begin{pmatrix}
      {\t P}^*&{\t Q}^*\\ 0&{\t S}^*
    \end{pmatrix}
    \begin{pmatrix}
      {\t P}&0\\ {\t Q}&{\t S}
    \end{pmatrix}
  \end{equation*}
  and $\t P$ chosen so that $S(S(1);0) = {\t P}^*\t P$.  Note that
  \begin{equation*}
    L =
    \begin{pmatrix}
      {\t S}^* & {\t T}^* \\ 0 & U^*
    \end{pmatrix}
    \begin{pmatrix}
      {\t S} & 0 \\ {\t T} & U
    \end{pmatrix}.
  \end{equation*}
  
  Since by assumption $\ran T \subseteq \clran U$, we have $S^*S =
  S(L;\{1\}) \geq {\t S}^*\t S$.  On the other hand, since
  \begin{equation*}
    S(1) \geq
    \begin{pmatrix}
      P^* & Q^* \\ 0 & S^*
    \end{pmatrix}
    \begin{pmatrix}
      P & 0 \\ Q & S
    \end{pmatrix},
  \end{equation*}
  we also have ${\t S}^*\t S \geq S^*S$.  Hence ${\t S}^* \t S= S^*S$.
  Thus $VS = \t S$ for some isometry $V : \clran S \to \clran \t S$.
  Since we have chosen $S(S(1);0) = {\t P}^*\t P$, by Lemma
  \ref{lemma1} $\ran {\t Q} \subseteq \clran \t S$.  Moreover,
  \begin{equation*}
    0 \leq
    \begin{pmatrix}
      {\t P}^* \t P + {\t Q}^* \t Q & {\t Q}^* S \\ S^* {\t Q} &
      S^*S
    \end{pmatrix}
    -
    \begin{pmatrix}
      P^*P + Q^*Q & Q^*S \\ S^*Q & S^*S
    \end{pmatrix}
  \end{equation*}
  and ${\t S}^*\t S \geq S^*S$ imply that $0={\t Q}^*\t S - Q^*S = (\t
  Q^* V - Q^*)S$.  As $\ran Q \subseteq \clran S$ it follows that $\t
  Q^* V = Q^*$.  Thus, in particular, $\t Q^* \t Q = Q^*Q$.  But then
  we obtain that
  \begin{equation}
    \label{eq:7}
    {\t P}^* \t P\geq P^*P.
  \end{equation}
  Observe that (\ref{eq:7}) will be true no matter what the original
  factorization of $M$ in (\ref{eq:3}) is as long as the range
  conditions in (\ref{eq:5}) are satisfied.
  
  Now instead consider the factorization $M = {G'}^*G'$, where
  \begin{equation*}
    G' =
    \begin{pmatrix}
      P'&0&0 \\ Q'&S&0 \\ R'&T&U
    \end{pmatrix}
  \end{equation*}
  with ${P'}^*P' = S(0)$.  Such a factorization is possible by Lemma
  \ref{lemma1}.  Since by assumption $\ran T \subseteq \clran U$, we
  have
  \begin{equation*}
    \clran
    \begin{pmatrix}
      S & 0 \\ T & U
    \end{pmatrix}
    \subseteq \clran S \oplus \clran U.
  \end{equation*}
  Also by Lemma \ref{lemma1} then,
  \begin{equation*}
    \ran
    \begin{pmatrix}
      {Q'} \\ {R'}
    \end{pmatrix}
    \subseteq \clran S \oplus \clran U,
  \end{equation*}
  and hence $\ran {Q'} \subseteq \clran S$.  So the conditions in
  (\ref{eq:5}) are satisfied for this factorization, and hence as
  noted above, we must have ${\t P}^*\t P \geq {P'}^*P'$.  But by
  definition of the Schur complement, ${P'}^*P' \geq {\t P}^*\t P$, so
  we have equality.  Consequently, (\ref{eq:4}) holds.
  
  Finally, using Lemma \ref{lemma1}, there is a factorization
  \begin{displaymath}
    L =
    \begin{pmatrix}
      S^* & T^* \\ 0 & U^*
    \end{pmatrix}
    \begin{pmatrix}
      S & 0 \\ T & U
    \end{pmatrix}
  \end{displaymath}
  where $\ran T \subseteq \clran U\subseteq {\mathcal H}_3$, so that
  $\clran \begin{pmatrix} S & T \\ 0 & U \end{pmatrix} = \clran S
  \oplus \clran U \subseteq {\mathcal H}_2 \oplus {\mathcal H}_3$.
  Again by Lemma \ref{lemma1}, there exists $P:{\mathcal H}_1 \to
  {\mathcal H}_1$ such that $P^*P = S(0)$ and (\ref{eq:3}) holds.
  Consequently $\ran\begin{pmatrix} Q \\ R \end{pmatrix} \subseteq 
  \clran S \oplus \clran U$, giving $\ran Q \subseteq \clran S$.
  
  It is now clear that the factorization in \eqref{eq:3} with these
  choices of $P,Q,R,S,T$ and $U$ satisfies the last statement of the
  theorem.
\end{proof}

\begin{corollary}
  \label{corollary1}
  Let $M\geq 0$ be an $n\times n$ operator matrix, $J\subseteq K
  \subseteq \{0\ldots n-1\}$.  Then
  \begin{equation}
    \label{eq:8}
    S(J) = S(S(K);J).
  \end{equation}
\end{corollary}

\begin{proof}
  Let $I_1= J$, $I_2 = K \setminus J$ and $I_3 = \{ 0, \ldots , n-1 \}
  \setminus K$.  Writing $M$ as a $3\times 3$ block matrix with
  respect to the partition
\begin{equation*}
\{ 0, \ldots , n-1
\} = I_1 \cup I_2 \cup I_3, 
\end{equation*}
 the corollary follows directly
Lemma \eqref{lemma2}.
\end{proof}

\begin{corollary}
  \label{corollary2}
  Given $M\geq 0$ an $n\times n$ operator matrix, there is a
  factorization $M = P^*P$ where
  \begin{displaymath}
    P =
    \begin{pmatrix}
      P_{00} & 0 & \cdots & \cdots & 0 \\
      P_{10} & P_{11} & 0 &\cdots& 0 \\
      \vdots & \vdots & \ddots & \ddots & \vdots \\
      \vdots & \vdots & \ddots & \ddots & 0 \\
      P_{n-1,0} & P_{n-1,1} & \cdots & \cdots & P_{n-1,n-1}
    \end{pmatrix},
  \end{displaymath}
  where $\clran P = \clran P_{00} \oplus \cdots \oplus \clran
  P_{n-1,n-1}$ and such that if $P_k$ is the truncation of $P$ to the
  upper left $(k+1)\times(k+1)$ corner, then $S(k) = P_k^*P_k$,
  $k=0,\ldots , n$.
\end{corollary}

The above result also appears in \cite{MR89a:47027}.

\begin{lemma} \label{lemma3}
  Let
  \begin{equation}\label{ppstar}
    \begin{pmatrix}
      P^* & Q^* \\ 0 & R^*
    \end{pmatrix}
    \begin{pmatrix}
      P & 0 \\ Q & R
    \end{pmatrix} =
    \begin{pmatrix}
      {\t P}^* & {\t Q}^* \\ 0 & {\t R}^*
    \end{pmatrix}
    \begin{pmatrix}
      \t P & 0 \\ \t Q & \t R
    \end{pmatrix},
  \end{equation}
  and suppose $\ran Q \subseteq
  \clran R$.  Then there is a unique isometry
  \begin{equation*}
    \begin{pmatrix}
      V_{11} & 0 \\ V_{21} & V_{22}
    \end{pmatrix}
  \end{equation*}
  acting on $\clran P \oplus \clran R$ so that
  \begin{equation} \label{io}
    \begin{pmatrix}
      \t P & 0 \\ \t Q & \t R
    \end{pmatrix} =  \begin{pmatrix} V_{11} & 0 \\ V_{21} & V_{22}
    \end{pmatrix} \begin{pmatrix}
      P & 0 \\ Q & R
    \end{pmatrix} .
  \end{equation}
\end{lemma}

\begin{proof} It is a standard result that $A^*A=B^*B$ if and only
  there exist an isometry $V :{\clran B} \to {\clran A}$ so that
  $VB=A$.  The operator $V$ is uniquely determined by setting $V(Bx) =
  Ax$ for every $Bx \in \ran B$, and extending $V$ to $\clran B$ by
  continuity.  Thus \eqref{ppstar} implies the existence of an
  isometry $V=(V_{ij})_{i,j=1}^2$ satisfying
  \begin{equation} \label{io2}
    \begin{pmatrix}
      \t P & 0 \\ \t Q & \t R
    \end{pmatrix} =  \begin{pmatrix} V_{11} & V_{12} \\ V_{21} & V_{22}
    \end{pmatrix} \begin{pmatrix}
      P & 0 \\ Q & R
    \end{pmatrix} .
  \end{equation}
  It remains to show that $V_{12}=0$.  Note that \eqref{io2} implies
  that $V_{22} R = \t R$.  Combining this with \eqref{ppstar} we get
  that
  \begin{equation*}
    R^* R = {\t R}^* \t R = { R}^* V_{22}^* V_{22} R ,
  \end{equation*}
  and thus
  \begin{equation}\label{iso}
    {R}^* (I_{\clran R} - V_{22}^* V_{22} ) R = 0.
  \end{equation}
  As $\ran Q \subseteq \clran R$ we have that
  \begin{equation*} 
    {\clran } \begin{pmatrix}
      P & 0 \\ Q & R
    \end{pmatrix} = {\clran P} \oplus {\clran R} .
  \end{equation*}
  Thus $V_{22}$ and $V_{12}$ act on ${\clran R}$.  From \eqref{iso} we
  now obtain that $V_{22}$ is an isometry on ${\clran R}$.  But then,
  since $V$ is an isometry, we must have that $V_{12} = 0$.
\end{proof}

In the following lemma we consider a positive operator on ${\mathcal
  H}_1 \oplus {\mathcal H}_2 \oplus {\mathcal H}_3 \oplus {\mathcal
  H}_4$, the ${\mathcal H}_k$'s Hilbert spaces.
\begin{lemma}\label{lemma4}
  Let $A=(A_{ij})_{i,j=1}^4 \ge 0$.  Then we have that
  \begin{equation}\label{zero}
    [S(2)]_{21}=0 \end{equation}
if and only if
\begin{equation}\label{schurcombo}
  S(2) = S(1) + S(\{0,2\}) - S(0).
\end{equation}
\end{lemma}

\begin{proof} The direction \eqref{schurcombo} $\Rightarrow$
  \eqref{zero} is trivial.
  
  By Corollary \ref{corollary2}, there is a lower triangular $3\times
  3$ operator matrix
  \begin{displaymath}
    P =
    \begin{pmatrix}
      P_{00} & 0 & 0 \\
      P_{10} & P_{11} & 0 \\
      P_{20} & P_{21} & P_{22} \\
    \end{pmatrix}
  \end{displaymath}
  such that $S(2) = P^*P$, $S(1) = \begin{pmatrix} P_{00}^* & P_{10}^* \\
    0 & P_{11}^* \end{pmatrix} \begin{pmatrix} P_{00} & 0 \\ P_{10} &
    P_{11} \end{pmatrix}$, and $S(0) = P_{00}^* P_{00}$.  Also, $\ran
  P_{10} \subseteq \clran P_{11}$ and $\ran P_{20}, \ran P_{21}
  \subseteq \clran P_{22}$.  Thus $[S(2)]_{21}=0 $ is equivalent to
  $P_{21} = 0$.  Interchanging the order of rows $1$ and $2$ and
  columns $1$ and $2$, we have
  \begin{displaymath}
    S(\{0,2,1\}) =
    \begin{pmatrix}
      P_{00}^* & P_{20}^* & P_{10}^* \\
      0 & P_{22}^* & 0 \\
      0 & 0 & P_{11}^* \\
    \end{pmatrix}
    \begin{pmatrix}
      P_{00} & 0 & 0 \\
      P_{20} & P_{22} & 0 \\
      P_{10} & 0 & P_{11} \\
    \end{pmatrix}
  \end{displaymath}
  Since $\ran \begin{pmatrix}P_{10} & 0 \end{pmatrix} \subseteq \clran
  P_{11}$, by Lemma \ref{lemma1},
  \begin{displaymath}
    S(\{0,2\}) =
    \begin{pmatrix}
      P_{00}^* & P_{20}^* \\
      0 & P_{22}^* \\
    \end{pmatrix}
    \begin{pmatrix}
      P_{00} & 0 \\
      P_{20} & P_{22} \\
    \end{pmatrix}.
  \end{displaymath}
  A direct calculation verifies the equality in \eqref{schurcombo}.
\end{proof}

By relabelling and grouping as we did in the proof of Corollary
\ref{corollary1}, we obtain the following.

\begin{corollary}
  \label{corollary3}
  Suppose $M \geq 0$ is an $n\times n$ operator matrix, $K\cup J = N
  \subseteq \{ 0, \ldots n-1\}$.  Then
  \begin{equation}
    \label{eq:9}
    S(N) = S(K) + S(J) - S(K\cap J)
  \end{equation}
  if and only if
  \begin{equation}
    \label{eq:10}
    \begin{split}
      [S(N)]_{k,j} = 0, \qquad(k,j) & \in (N\times N) \setminus
      ((K\times K) \cup (J\times J))\\
      & = [(K\setminus (K\cap J))\times (J\setminus (K\cap J))] \cup
      [(J\setminus (K\cap J)\times (K\setminus (K\cap J))]
    \end{split}
  \end{equation}
\end{corollary}

\section{One variable outer and inner-outer factorization}
\setcounter{equation}{0}

In this section we will provide new proofs for several one variable
factorization results.  These proofs are based on the properties of
Schur complements.

Given a Hilbert space $\sp H$ let $\HardyD H$ denote the Hardy space
of $\sp H$-valued functions analytic in the unit disk with square
integrable boundary values.  These functions will be identified with
their boundary values whenever convenient.  Given a pair of Hilbert
space $\sp H$, $\sp K$, let $\L HK$ stand for the Banach space of
bounded operators acting ${\sp H} \to {\sp K}$.  We will write ${\bf
  L} ({\sp H})$ instead of ${\bf L}({\sp H},{\sp H})$.  As usual,
$\HinfDL HK$ stands for the set of all bounded holomorphic $\L HK
$-valued functions on $\mathbb D$.  With the operator valued function
$F\in\HinfDL HK$, we associate the operator $M_F:\HardyD H\to \HardyD
K$ of multiplication by $F$; that is, $M_Fg(z)=F(z)g(z)$.  The
function $F$ is called {\it outer\/} if the corresponding
multiplication operator $M_F$ has dense range in $\HardyD M$ for some
subspace $\sp M$ of $\sp K$, and this reduces to the usual definition
when $\mathcal H$ and $\mathcal K$ are $\mathbb C$.  For $Q \in
L^\infty_{\Ls H}({\mathbb T})$, we consider the Toeplitz operator
$T_Q:\HardyD H\to \HardyD H$ defined via $T_Q f = \Pi_+ (Qf)$, where
$\Pi_+$ the projection is from $L^2_{\mathcal H} ({\mathbb T})$ onto
$\HardyD H$.  We shall often represent $T_Q$ via the Toeplitz operator
matrix
\begin{equation} \label{toepmatrix}
  T_Q \equiv
  \begin{pmatrix} Q_0 & Q_{-1} & \cdots \cr Q_{1} & \ddots &
    \ddots \cr \vdots & \ddots & & \end{pmatrix},
\end{equation}
where we make the obvious identification of $f(z) = \sum_0^\infty f_k
z^k \in \HardyD H$ with $\mathrm{col}(f_k)_0^\infty\in
\ell^2_{\mathcal H} ({\mathbb N}_0)$, where $f_j \in {\mathcal H}$ and
$\| f\| := \sqrt{ \sum_{j=0}^\infty \| f_j \|^2} < \infty $.  We view
the matrix as an operator matrix with rows and columns indexed by
${\mathbb N}_0 = \{ 0,1, \ldots \}$.  In addition, we shall often use
the identification
\begin{equation} 
  \label{ident} 
  T_Q = \begin{pmatrix} Q_0 & {\rm row} (Q_{-i})_{i\ge 1} 
    \cr {\rm col} (Q_i)_{i\ge 1}& T_Q \end{pmatrix}.
\end{equation}
In other words, the operator $L: z \HardyD H \to z \HardyD H$ defined
by $(Lf)(z) = z \Pi_+(Qz^{-1}f)$ will at times be identified with
$T_Q$.  Following the notation from the previous section, for $\Lambda
\subset {\mathbb N}_0$ we let $S(T_Q; \Lambda )$ (or $S(\Lambda )$
when no confusion is possible) denote the Schur complement of $T_Q$
supported on rows and columns indexed by $\Lambda$.  In addition,
$S(k)$ is a shorthand for $S(\{ 0,\ldots , k \})$.
  
We first address the effect the Toeplitz structure has on the Schur
complements.

\begin{proposition} \label{scs} Consider the positive semidefinite Toeplitz
  operator $T_Q = (Q_{i-j})_{i,j=0}^\infty$ acting on $\ell^2_{\mathcal
    H} ({\mathbb N}_0)$.  Then the Schur complements $S(m)$ of $T_Q$
  satisfy the recurrence relation
  \begin{equation} \label{form}
    S(m) = \begin{pmatrix} A & B^* \cr B & S(m-1) \end{pmatrix} ,
  \end{equation}
  for appropriate choice of $A: {\mathcal H} \to {\mathcal H}$ and $B:
  {\mathcal H} \to {\mathcal H}^m$.  When $Q_j = 0, j\ge m+1$, then
  $A=Q_0$ and $B={\rm col} (Q_i)_{i=1}^m$.
\end{proposition}

\begin{proof}
  By the definition of Schur complement
  \begin{equation} \label{form2}
    T_Q - \begin{pmatrix} S(m) & 0 \cr 0 & 0 \end{pmatrix} \ge 0 .
  \end{equation}
  Let us write 
  \begin{equation*}
    S(m) = \begin{pmatrix} A & B^* \cr B & C \end{pmatrix} :
    {\mathcal H} \oplus {\mathcal H}^m 
    \to {\mathcal H} \oplus {\mathcal H}^m.  
  \end{equation*}
  Leaving out row and column 0 in \eqref{form2} yields
  \begin{equation*}
    T_Q - \begin{pmatrix} C & 0 \cr 0 & 0 \end{pmatrix} \ge 0 , 
  \end{equation*}
  where we used identification \eqref{ident}.  This shows that $C \le
  S(m-1)$.  On the other hand, leaving out row and columns $1,\ldots ,
  m$ in \eqref{form2} yields
  \begin{equation*}
    \begin{pmatrix} Q_0 - A & {\rm row} (Q_j^*)_{j\ge m+1} \cr
      {\rm  col} (Q_j)_{j\ge m+1} & T_Q \end{pmatrix} \ge 0 .  
  \end{equation*}
  Hence
  \begin{equation*}
    A \le S(\begin{pmatrix} Q_0 & {\rm row} (Q_j^*)_{j\ge m+1} \cr
      {\rm  col} (Q_j)_{j\ge m+1} & T_Q \end{pmatrix} ; 0 ) =: \tilde A . 
  \end{equation*}
  
  Note that when $Q_j=0$, $j \ge m+1$, we have that $\tilde A = Q_0$.
  Consider now the operator matrix
  \begin{equation}\label{atilde}
    \begin{pmatrix} Q_0 - \tilde A & X & {\rm row} (Q_j^*)_{j\ge m+1} \cr X^* &
      (Q_{i-j})_{i,j=1}^m - S(m-1) & (Q_{i-j})_{i=1, j=m+1}^{m+1,\ \infty} \cr
      {\rm  col} (Q_j)_{j\ge m+1}& (Q_{i-j})_{i=m+1, j=1}^{\infty,\ m
      } & T_Q 
    \end{pmatrix}.  
  \end{equation}
  The existence of an operator $X$ making this into a positive
  semidefinite matrix is a variant of a standard operator matrix
  completion problem, and by \cite{MR89k:47010} (see also, e.g.,
  Theorem XVI.3.1 in \cite{MR92k:47033} or \cite{MR94j:47024}), there
  is always such an $X$.  Note that when $\tilde A= Q_0$ we have necessarily
  that $X = 0$.  As \eqref{atilde} is positive semidefinite we obtain that
  \begin{equation*}
    \begin{pmatrix} \tilde A & {\rm row} (Q_j^*)_{j=1}^m - X \cr 
      {\rm col} (Q_j)_{j=1}^m -X^*
      & S(m-1) \end{pmatrix} \le S(m) = 
    \begin{pmatrix} A & B^* \cr B& C \end{pmatrix} . 
  \end{equation*}
  This implies that $\tilde A \le A$ and $S(m-1) \le C$.  As we also
  had that $A \le \tilde A$ and $C \le S(m-1)$, the equalities $A = \tilde A$ and $C = S(m-1)$ 
follow.  This yields \eqref{form}.  Moreover, when
  $Q_j=0$ for $j\ge m+1$, we have that $\tilde A=Q_0$ and $X=0$, and
  thus $B={\rm col} (Q_i)_{i=1}^m$.
\end{proof}

{\bf Remark.} Note that the proof shows that the operator $A$ in
\eqref{form} is given by
\begin{equation*}
 A =  S\left(\begin{pmatrix} Q_0 & {\rm row} (Q_j^*)_{j\ge m+1} \cr
{\rm  col} (Q_j)_{j\ge m+1} & T_Q \end{pmatrix} ; 0 \right) . 
\end{equation*}

Because of the inheritance principle observed in Proposition
\ref{scs}, the Schur complements of a Toeplitz operator allow a
stationary $UL$ Cholesky decomposition.

\begin{corollary} \label{scs3} Consider the positive semidefinite Toeplitz
  operator $T_Q = (Q_{i-j})_{i,j=0}^\infty$ acting on $\ell^2_{\mathcal
    H}({\mathbb N}_0)$.  Then there exist operators $F_0 , F_1 ,
  \ldots $ with $F_i : {\mathcal H} \to \clran F_0 \subseteq {\mathcal
    H}$ so that the Schur complements $S(m)$ of $T_Q$ satisfy
  \begin{equation} \label{stat}
    S(m) = \begin{pmatrix} F_0^* & \cdots & F_m^* \cr & \ddots &
      \vdots \cr & & F_0^*
    \end{pmatrix}
    \begin{pmatrix} F_0 & & \cr \vdots & \ddots & \cr F_m & \cdots &
      F_0 
    \end{pmatrix}
    , \qquad m\ge 0.
  \end{equation}
\end{corollary}

\begin{proof}
  We prove this by induction.  When $m=0$ we may for instance choose
  $F_0 = (S(0))^{1/2}$.  It follows from Proposition \ref{scs} that
  $(S(m))_{m,m} = (S(m-1))_{m-1,m-1} = F_0^* F_0$, where in the last
  step we used the induction hypothesis.  By Corollary
  \ref{corollary1} we have that $S(m-1) = S(S(m);m-1)$, and thus by
  Lemma \ref{lemma1} with
  \begin{equation*}
    P = 
    \begin{pmatrix} F_0 & & \cr \vdots & \ddots & \cr F_{m-1} & 
      \cdots & F_0 \end{pmatrix} 
    ,\qquad R=F_0, 
  \end{equation*}
  there exist $\begin{pmatrix} G_m & \cdots & G_1 \end{pmatrix}$ so
  that
  \begin{equation}\label{form3}
    S(m) = \begin{pmatrix} F_0^* & \cdots & F_{m-1}^* & G_m^* \cr 
      & \ddots & \vdots & \vdots \cr
      & & F_0^* & G_1^* \cr & & & F_0^* \end{pmatrix}
    \begin{pmatrix} F_0 & & & \cr \vdots & \ddots &  &  \cr
      F_{m-1} & \cdots & F_0 & \cr G_{m} & \cdots & G_1 & F_0 \end{pmatrix} ,
  \end{equation}
  and $\ran \begin{pmatrix} G_m & \cdots & G_1 \end{pmatrix} \subseteq
  \clran F_0$.  Comparing \eqref{form3} with \eqref{form} along with
  the induction hypothesis yields
  \begin{equation*}
    \begin{split}
      &\begin{pmatrix} F_0^* & \cdots & F_{m-2}^* & G_{m-1}^* \cr 
        & \ddots & \vdots & \vdots \cr
      & & F_0^* & G_1^* \cr & & & F_0^* \end{pmatrix}
    \begin{pmatrix} F_0 & & & \cr \vdots & \ddots &  &  \cr
      F_{m-2} & \cdots & F_0 & \cr G_{m-1} & \cdots & G_1 & F_0
    \end{pmatrix} 
    = S(m-1) = \\
    &\qquad\qquad = \begin{pmatrix} F_0^* & \cdots & F_{m-2}^* &
      F_{m-1}^* \cr 
      & \ddots & \vdots & \vdots \cr
      & & F_0^* & F_1^* \cr & & & F_0^* \end{pmatrix}
    \begin{pmatrix} F_0 & & & \cr \vdots & \ddots &  &  \cr
      F_{m-2} & \cdots & F_0 & \cr F_{m-1} & \cdots & F_1 & F_0 \end{pmatrix} ,
    \end{split}
  \end{equation*}
  and thus
  \begin{equation*}
    F_0^*  \begin{pmatrix} G_{m-1} & \cdots & G_1 \end{pmatrix} =
    F_0^*  \begin{pmatrix} F_{m-1} & \cdots & F_1 \end{pmatrix} . 
  \end{equation*}
  As $\ran \begin{pmatrix} G_{m-1} & \cdots & G_1 \end{pmatrix}
  \subseteq \clran F_0 $ and $\ran \begin{pmatrix} F_{m-1} & \cdots &
    F_1 \end{pmatrix} \subseteq \clran F_0 $, it follows that $G_j =
  F_j$, $j=1,\ldots , m-1$.  By setting $F_m := G_m$, we obtain the
  result.
\end{proof}

Before we come to our main results, let us develop some equivalent
statements for outerness that follow directly from the Schur
complement results.  Analogously to \eqref{ident}, we shall use the
identification
\begin{equation}\label{mf} 
  T_F = \begin{pmatrix} F_0 & 0 \cr {\rm col} (F_j)_{j\ge 1} & T_F
  \end{pmatrix} . 
\end{equation} 

\begin{theorem} \label{eqs} Let $F\in\HinfDL HK$.  Denote the Taylor
  coefficients of $F$ by $F_j$, $j\ge 0$.  The following are
  equivalent:
  \begin{enumerate}[(i)]
  \item $F$ is outer;
  \item $\clran M_F = H^2_{\clran F_0} ({\mathbb D})$;
  \item $\ran \ {\rm col} (F_j)_{j\ge 1} \subset \clran T_F$;
  \item $S(T_F^* T_F; 0) = F_0^* F_0$
  \item For some $k \in {\mathbb N}_0$ we have that
    \begin{equation} \label{sc} 
      S(T_F^* T_F ; k) = \begin{pmatrix} F_0^* & \cdots & F_k^* \cr 
        & \ddots & \vdots \cr & & F_0^* \end{pmatrix}  
      \begin{pmatrix} F_0 &  &  \cr \vdots
        & \ddots & \cr F_k & \cdots & F_0 \end{pmatrix} ; 
    \end{equation}
  \item For all $k \in {\mathbb N}_0$ equality \eqref{sc} holds;
\end{enumerate}
\end{theorem}

\begin{proof} Clearly (ii) implies (i).  For the implication (i)
  $\Rightarrow$ (ii), observe that if $\clran M_F = H^2_{\sp M}
  ({\mathbb D})$, then $P_0 (\clran T_F ) = {\sp M}$ where $P_0$ is
  the projection $F \to F_0$.  But when $h \in H^2_{\sp M} ({\mathbb
    D})$ we have that $P_0 (Fh) = F_0 h(0) \in {\ran F_0}$.  Moreover,
  since we may let $h(0)$ range over all elements in ${\sp H}$, we
  obtain that ${\sp M}= \clran F_0$.
  
  For (ii) $\Rightarrow$ (iii), note that given (ii) we get that $\ran
  {\rm col} (F_j)_{j\ge 1} \subset (I-P_0) \ell^2_{\clran
    F_0}({\mathbb N}_0) \approx H^2_{\clran F_0} ({\mathbb D}) =
  \clran T_F $.  We used here the identification a multiplication
  operator and its corresponding Toeplitz operator.
  
  Next consider (iii) $\Rightarrow$ (ii).  If $h_0 \in {\sp H}$ we
  have that ${\rm col} (F_j h_0)_{j\ge 1} \in \clran T_F$.  Thus there
  exist $g_i$ so that $\lim_{i\to \infty} T_F g_i = {\rm col} (F_j
  h_0)_{j\ge 1}$.  But then
  \begin{equation*}
    \begin{pmatrix} F_0 & 0 \cr {\rm col}
      (F_j)_{j\ge 1} & T_F
    \end{pmatrix} \begin{pmatrix} h_0 \cr -g_i \end{pmatrix}
    \to \begin{pmatrix} F_0 h_0 \cr 0 \end{pmatrix}.
  \end{equation*}
  Thus $F_0 h_0 \in \clran M_F$.  As $\clran M_F$ is closed under
  multiplication with $z$, we get that $H^2_{\clran F_0} ({\mathbb D})
  \subseteq {\clran} M_F$.  The inclusion $\clran M_F \subseteq
  H^2_{\clran F_0} ({\mathbb D})$ follows as (iii) implies $\ran F_j
  \subseteq \clran \begin{pmatrix} F_{j-1} & \cdots & F_0
  \end{pmatrix}$, $j \ge 1$, which in turn implies $\ran F_j \subseteq
  \clran F_0$.
  
  For the implication (iii) $\Rightarrow$ (vi), note that (iii)
  implies that $\ran {\rm col} (F_j)_{j\ge k} \subset \clran T_F$ for
  all $k\ge 1$.  But then (vi) follows immediately from Lemma
  \ref{lemma1}.
  
  The implications (vi) $\Rightarrow$ (iv) $\Rightarrow$ (v) are
  trivial.
  
  For (v) to (iii) use Lemma \ref{lemma1}.
\end{proof}

We are now ready to give a simple proof for the operator valued
Fej\'er-Riesz theorem.  The original proof is due to Rosenblum
\cite{MR37:3378}.

\begin{theorem} \label{RF}
  (\cite{MR37:3378}) Let $Q_j: {\sp H} \to {\sp H}$, $j=-m,\ldots ,
  m$, be Hilbert space operators so that $Q(z) := \sum_{j=-n}^n Q_j
  z^j \ge 0$, $z\in\mathbb T$.  Then there exists an outer operator
  polynomial $P(z) = \sum_{j=0}^m P_j z^j$ with $P_j \in \Ls H$,
  $j=0,\ldots , m$, so that $Q(z) = P(z)^* P(z)$, $z \in \mathbb T$.
\end{theorem}

\begin{proof}
  Let
  \begin{equation*}
    Y= (Y_{ij})_{i,j=0}^m :=S(m) - S(m-1),
  \end{equation*}
  where $S(m-1)$ is viewed as an operator on ${\mathcal H}^{m+1}$
  (with last row and column equal to 0).  By Corollary \ref{scs3} we
  have that there exist operators $P_i: {\mathcal H} \to {\mathcal H}$
  with $\ran P_i \subseteq \clran P_0$ so that
  \begin{equation*}
    Y = \begin{pmatrix} P_m^* \cr \vdots \cr P_0^* \end{pmatrix}
    \begin{pmatrix} P_m & \cdots & P_0 \end{pmatrix} .
  \end{equation*}
  Put $P(z) =\sum_{j=0}^m P_j z^j$, $Z_m = {\begin{pmatrix} z^m &
      \cdots & 1 \end{pmatrix}}^T$.  Then, since in Proposition
  \ref{scs} we have that $A=Q_0$ and $B={\rm col}(Q_i)_{i=1}^m$, we
  get that
  \begin{equation*}
    \begin{split}
      P(z)^* P(z) &= Z_m^* Y Z_m \\
      &= Q(z) + Z_{m-1}^* S(m-1) Z_{m-1} - \overline{z}Z^*_{m-1} S(m-1)
      zZ_{m-1},
    \end{split}
  \end{equation*}
  where the last two terms cancel when $z\in\mathbb T$.
  
  Finally, in order to see that $P(z)$ is outer, use the equivalence
  (i) $\Leftrightarrow$ (iv) in Theorem \ref{eqs} and the fact that
  \begin{equation*}
    S(T_P^* T_P; 0)= S(T_Q; 0)= S(S(m);0) = P_0^* P_0.
  \end{equation*}
\end{proof}

Next we will show how Lemma \ref{lemma3} leads to the existence of
inner-outer factorizations for operator valued polynomials.  Recall
that $A \in H^\infty_{\L HK } ({\mathbb T})$ is {\it inner} if the
multiplication operator $M_A:H^2_{\mathcal H} ({\mathbb T}) \to
H^2_{\mathcal K} ({\mathbb T}) $ with symbol $A$ is a partial
isometry.

\begin{theorem}[Existence of inner-outer factorization] \label{exists}
  Let $A\in H^\infty_{\L HK} ({\mathbb D})$.  Then there exists an
  outer function $F$ and an inner function $V$, so that $A=VF$.
\end{theorem}

\begin{proof}
  Consider the Toeplitz operator $T_Q := T_A^* T_A$.  By Corollary
  \ref{scs3} there exist $F_j :{\mathcal H} \to \clran F_0 \subseteq
  {\mathcal H}$ so that
  \begin{equation} \label{stat2}
    S(m) = \begin{pmatrix} F_0^* & \cdots & F_m^* \cr 
      & \ddots & \vdots \cr & & F_0^* \end{pmatrix}
    \begin{pmatrix} F_0 & & \cr \vdots & \ddots & \cr 
      F_m & \cdots & F_0 \end{pmatrix} =:
    {\mathcal F}(m)^*{\mathcal F}(m) ,\qquad m\ge 0.
  \end{equation}
  Note also that
  \begin{equation} \label{stat3}
    {\mathcal A}(m)^*{\mathcal A}(m) := 
    \begin{pmatrix} A_0^* & \cdots & A_m^*
      \cr & \ddots & \vdots \cr & & A_0^* \end{pmatrix}
    \begin{pmatrix} A_0 & & \cr \vdots & \ddots & \cr 
      A_m & \cdots & A_0 \end{pmatrix} \le S(m).
  \end{equation}
  Consider the sequence of operators
  \begin{equation}\label{star}
    \begin{pmatrix} {\mathcal F}(m) & 0 \cr 0 & 0 \end{pmatrix}
  \end{equation}
  acting on $\ell^2_{\mathcal H}({\mathbb N}_0)$.  As $\|S(m) \| \le\|
  T_A^* T_A\|$, it follows that \eqref{star} is a bounded sequence of
  operators, and therefore has a subsequence that converges to $T_F$,
  say, in the weak-$*$ topology.  But then we must have that
  \begin{equation*}
    T_F = \begin{pmatrix} F_0 & & & & \cr F_1 & F_0 & & & \cr 
      \vdots & \ddots & \ddots & &
    \end{pmatrix}. 
  \end{equation*}
  
  Also, ${\mathcal A}(m)$ converges to $T_A$ in the weak-$*$ topology.
  But now \eqref{stat3} and $\begin{pmatrix} S(m) & 0 \cr 0 & 0
  \end{pmatrix} \le T_A^* T_A $ yield that $T_A^* T_A \le T_F^* T_F
  \le T_A^* T_A$.  Thus $T_A^* T_A = T_F^* T_F$, or equivalently, $
  A(z)^* A(z) = F(z)^* F(z)$ a.e on ${\mathbb T}$, where $F(z) = F_0 +
  z F_1 + \ldots $.  As $S(T_F^* T_F; 0 ) = F_0^* F_0$ it follows by
  Theorem \ref{eqs} that $F$ is outer.

  Next, notice that we may write
  \begin{equation*}
    T_F = \begin{pmatrix} F_0 & 0 \cr {\rm col} (F_j)_{j\ge 1} & T_F \end{pmatrix},
    \qquad
    T_A = \begin{pmatrix} A_0 & 0 \cr {\rm col} (A_j)_{j\ge 1} & T_A \end{pmatrix} . 
  \end{equation*}
  
  Moreover, since $F$ is outer we have that $\ran {\rm col}
  (F_j)_{j\ge 1} \subset \clran T_F$.  By Lemma \ref{lemma3} there
  exists a unique isometry
  \begin{equation*}
    \tilde  V = \begin{pmatrix}
      V_{11} & 0 \\ V_{21} & V_{22}
    \end{pmatrix}
  \end{equation*}
  acting on $\clran F_0 \oplus \clran T_F$ so that
  \begin{equation*}
    \begin{pmatrix}
      A_0 & 0 \\ {\rm col} (A_j)_{j\ge 1} & T_A
    \end{pmatrix} =  \begin{pmatrix} V_{11} & 0 \\ V_{21} & V_{22}
    \end{pmatrix} \begin{pmatrix}
      F_0 & 0 \\ {\rm col} (F_j)_{j\ge 1} & T_F
    \end{pmatrix} .
  \end{equation*}
  Since $V_{22}$ is an isometry and satisfies $T_A = V_{22} T_F$ we
  obtain by the uniqueness statement in Lemma \ref{lemma3} that
  $\tilde V=V_{22}$.  But that implies that $\tilde V$ must be of the
  form $ \tilde V= (V_{i-j})_{i,j\ge 0}$ with $V_k = 0$ for $k <0$.
  Thus $\tilde V = T_V$ and $V(z) = V_0 + z V_1 + z^2 V_2 + \ldots $
  is inner.
\end{proof}

Next we provide new proofs to some more of the various equivalent
characterizations that exist for outer functions (see
\cite{MR87e:47001}), and obtain a few new ones as well.

\begin{theorem}\label{eqs2} Let $F\in\HinfDL HK$.  Denote the Taylor
  coefficients of $F$ by $F_j$, $j\ge 0$.  The following are
  equivalent:
  \begin{enumerate}[(i)]
  \item $F$ is outer;
  \item For any $z\in\mathbb D$ and $G \in\HinfDL LK$ with $G^*G = F^*
    F$ a.e.~on $\mathbb T$,
    \begin{equation*}
      G(z)^* G(z) \le F(z)^* F(z) , z \in {\mathbb D} ; 
    \end{equation*}
  \item There exists $z_0 \in {\mathbb D}$ such that whenever $G
    \in\HinfDL LK$ and $G^*G = F^* F$ a.e.~on $\mathbb T$,
    \begin{equation*}
      G(z_0)^* G(z_0) \le F(z_0)^* F(z_0) ; 
    \end{equation*}
  \item $G \in\HinfDL LK$ and $G^*G = F^* F$ a.e.~on $\mathbb T$
    implies
    \begin{equation*}
      G_0^* G_0 \le F_0^* F_0 ; 
    \end{equation*}
  \item For some $k \in {\mathbb N}_0$ we have that $G \in\HinfDL LK$
    and $G^*G = F^* F$ a.e.~on $\mathbb T$ implies $\sum_{i=0}^l G_i^*
    G_i \le \sum_{i=0}^l F_i^* F_i$, $l=0,\ldots , k$, where $G_i$ are
    the Taylor coefficients of $G$;
  \item For all $k \in {\mathbb N}_0$ we have that $G \in\HinfDL LK$
    and $G^*G = F^* F$ a.e.~on $\mathbb T$ implies $\sum_{i=0}^k G_i^*
    G_i \le \sum_{i=0}^k F_i^* F_i$, where $G_i$ are the Taylor
    coefficients of $G$.
  \item For some $k \in {\mathbb N}_0$ we have that $G \in\HinfDL LK$
    and $G^*G = F^* F$ a.e.~on $\mathbb T$ implies
    \begin{equation}\label{mg1} 
      \begin{pmatrix} F_0^* & \cdots & F_k^* \cr & \ddots & \vdots
        \cr & & F_0^* \end{pmatrix}  \begin{pmatrix} F_0 &  &  \cr \vdots
        & \ddots & \cr F_k & \cdots & F_0 \end{pmatrix}\ge
      \begin{pmatrix} G_0^* & \cdots & G_k^* \cr & \ddots & \vdots
        \cr & & G_0^* \end{pmatrix}  \begin{pmatrix} G_0 &  &  \cr \vdots
        & \ddots & \cr G_k & \cdots & G_0 \end{pmatrix} , 
    \end{equation}
    where $G_i$ are the Taylor coefficients of $G$;
  \item For all $k \in {\mathbb N}_0$ we have that $G \in\HinfDL LK$
    and $G^*G = F^* F$ a.e.~on $\mathbb T$ implies \eqref{mg1}.
  \end{enumerate}
\end{theorem}

\begin{proof}
  For any $G \in\HinfDL LK$ we have that
  \begin{equation}\label{mg} 
    S(T_G^* T_G ; k) \ge
    \begin{pmatrix} G_0^* & \cdots & G_k^* \cr & \ddots & \vdots
      \cr & & G_0^* \end{pmatrix} \begin{pmatrix} G_0 & & \cr \vdots &
      \ddots & \cr G_k & \cdots & G_0 
    \end{pmatrix} . 
  \end{equation}
  Combining this observation with Theorem \ref{eqs} (v) and the fact
  that $T_F^*T_F = T_G^*T_G$, we immediately obtain the implication
  (i) $\Rightarrow$ (viii).
  
  The implications (viii) $\Rightarrow$ (iv) $\Rightarrow$ (vii)
  $\Rightarrow$ (v) $\Rightarrow$ (iv), (viii) $\Rightarrow$ (vi)
  $\Rightarrow$ (v) and (ii) $\Rightarrow$ (iv) $\Rightarrow$ (iii)
  are trivial.
  
  For (iv) $\Rightarrow$ (i), let $F=V\tilde F$ be an inner-outer
  factorization of $F$.  Then, by Theorem \ref{eqs} we have that
  \begin{equation*}
    F_0^* F_0 \le  S(T_F^* T_F; 0 ) = S(T_{\tilde F}^* T_{\tilde F};
    0) 
    = \tilde F_0^* \tilde F_0 . 
  \end{equation*}
  On the other hand, by (iv) $\tilde F_0^* \tilde F_0 \le F_0^* F_0$,
  and thus equality $F_0^* F_0 = S(T_F^* T_F; 0 )$ holds.  Again
  applying Theorem \ref{eqs}, gives that $F$ is outer.
  
  For the implication (i) $\Rightarrow$ (ii) fix $z\in\mathbb D$ and
  introduce the Blaschke factor $b_{z}(w) = -{w-z \over 1 - \bar{z} w
  }$, $w\in \mathbb D$.  Then $F$ is outer if and only if $F\circ b_z$
  is (use that the composition operators $g \to g \circ b_z$ and $g \to g\circ b_z^{-1}$
are bounded operators on $H_{\mathcal H}^2 
({\mathbb D})$ \cite[Theorem 3.6]{CM}).  Moreover $F(w)^*F(w) = G(w)^*G(w)$ a.e.~on ${\mathbb T}$ if and
  only if $F(b_z(w))^*F(b_z(w)) = G(b_z(w))^*G(b_z(w))$ a.e.~on
  ${\mathbb T}$.  Since $F \circ b_z$ is outer, by (i) $\Rightarrow$
  (iv), we have that
  \begin{equation*}
    F(b_z(0))^* F(b_z(0)) \ge \tilde G_0^* \tilde G_0 
  \end{equation*}
  for any $\tilde G$ such that $\tilde G^* \tilde G = (F \circ b_z)^*
  (F \circ b_z)$ a.e.~on ${\mathbb T}$.  Putting now $G = \tilde G
  \circ b_z^{-1}$ gives that
  \begin{equation*}
    F(z)^* F(z) \ge G(z)^* G(z)
  \end{equation*}
  for any $G$ with $G^*G = F^*F$ a.e.~on ${\mathbb T}$.  Since $z\in
  {\mathbb D}$ was arbitrary, the result follows.
  
  As (iv) $\Rightarrow$ (i) holds, it follows that if $F$ satisfies
  (iii), then $F\circ b_{z_0}$ is outer.  But then $F$ is outer as
  well.  This proves (iii) $\Rightarrow$ (i).
\end{proof}

\section{Multivariate outer polynomials}
\setcounter{equation}{0}

With the ideas from the previous section we now present a multivariate
operator-valued version of the Fej\'er-Riesz lemma.  As mere positive
semidefiniteness on the $d$-torus does not suffice, an additional
condition is required for $Q$ to allow an ``outer'' factorization.
This additional condition on $Q$ is given in terms of Schur
complements of $T_Q$, the Toeplitz operator on ${H^2_{\sp H}({\mathbb
    D}^d)}$ with symbol $Q$.

In order to state the result precisely we need some additional
notation.  For $z=(z_1,\ldots , z_d) \in {\mathbb T}^d$ and $k =
(k_1,\ldots , k_d ) \in {\mathbb Z}^d$ define $z^k := z_1^{k_1} \cdots
z_d^{k_d}$.  In this case $z^{*k} = \overline{z}^k = z^{-k}$.  We
write $0$ for $(0,\ldots,0)$.  For set $A, B \subseteq {\mathbb Z}^d$
we denote $A-B=\{ a-b \ : \ a \in A, b \in B \}$.  For matrices
labelled by elements of ${\mathbb Z}^d$ we fix the ordering as
lexicographical.  Since this is a total ordering, various results from
the first section on Schur complements readily translate to this
setting.  As before, we use the notation $S(T_Q; \Lambda )$ (or simply
$S(\Lambda )$ when no confusion is likely) to indicate a Schur
complement of $T_Q$ supported in rows and columns $\Lambda \subseteq
{\mathbb N}_0^d$.  In the same manner as when we labelled matrices
using elements of ${\mathbb N}_0$, we sometimes pad Schur complements
with zeros.  In this way for example, if $\Lambda_2 \subseteq
\Lambda_1$, then $S(\Lambda_1) - S(\Lambda_2)$ makes sense.  Finally,
we need the projections $\Pi_K$, $K\subseteq {\mathbb N}_0^d$, on
${H^2_{\sp M}({\mathbb D}^d)}$ defined by
\begin{equation*}
  \Pi_K \left(\sum_{k \in{\mathbb N}_0^d} h_k z^k \right) 
  = \sum_{k \in K} h_k z^k .
\end{equation*}

\begin{theorem} \label{multi}
  Let $K=\prod_{i=1}^d \{ 0,\ldots , n_i \}$ and let $Q_k: {\sp H} \to
  {\sp H}$, $k \in K - K$, be Hilbert space operators so that $Q(z) :=
  \sum_{k \in K-K} Q_k z^k \ge 0, z \in {\mathbb T}^d$.  Furthermore,
  let $n = (n_1 , \ldots , n_d )$ and $Z_K$ be the column matrix
  $(z^{n-k})_{k\in K}$.  The following are equivalent:
  \begin{enumerate}[(i)]
  \item there exists an operator polynomial $P(z) = \sum_{k\in K} P_k
    z^k$ with $P_j \in \Ls H$, $j \in K$, so that $Q(z) = P(z)^* P(z),
    z \in {\mathbb T}^d $, and
    \begin{equation}
      \label{inclusion1}
      \ran (\Pi_{{\mathbb N}_0^d \setminus K}\, T_P \, \Pi_{\{ n \} })
      \subseteq \clran (\Pi_{{\mathbb N}_0^d \setminus K} \, T_P\,
      \Pi_{{\mathbb N}_0^d \setminus K}) ,
    \end{equation}
    and
    \begin{equation}
      \label{inclusion2}
      \ran P_k \subseteq \clran P_0,\qquad  k\in K; \end{equation}
  \item The operator
    \begin{equation*}
      Y:= S(K ) - S( K \setminus \{ n \} )
    \end{equation*}
    satisfies
    \begin{equation}
      \label{condY}
      Z_K^* Y Z_K = Q(z),\qquad z \in {\mathbb
        T}^d . \end{equation}
  \end{enumerate}
\end{theorem}

\begin{proof} Suppose (ii) holds.  By Lemma \ref{lemma2} there exist
  $P_k \in \Ls H$, $k\in K$ such that with $P_K =
  \mathrm{row}(P_k)_{k\in K}$, $Y= P_K^*P_K$.  Defining $P(z) = P_KZ_K
  = \sum_{k\in K}P_k z^k$, we obtain from \eqref{condY} that $P(z)^*
  P(z) = Q(z), z \in {\mathbb T}^d $.  But then $T_Q = T_P^* T_P$.  View
  this factorization of $T_Q$ with respect to the decomposition
  \begin{equation}
    \label{eq:11}
    \ran \Pi_{K \setminus \{ n \} } \oplus \ran \Pi_{\{ n \} }
    \oplus \ran \Pi_{ {\mathbb N}_0 \setminus K },
  \end{equation}
  in which respect $T_P$ is a $3 \times 3$ lower triangular operator
  matrix.  We are now exactly in the situation of Lemma \ref{lemma2}
  with $T^* = \Pi_{{\mathbb N}_0^d \setminus K}\, T_P\, \Pi_{\{ n \}
  }, U^* = \Pi_{{\mathbb N}_0^d \setminus K}\, T_P\, \Pi_{{\mathbb
      N}_0^d \setminus K} , Q^* = \Pi_{\{ n \} }\, T_P\, \Pi_{K
    \setminus \{ n \} }$, and $ S^* =\Pi_{\{ n \} }\, T_P\,
  \Pi_{\{ n \} } = P_0 .  $ Since \eqref{eq:4} in Lemma \ref{lemma2}
  holds, we obtain \eqref{eq:5} of Lemma \ref{lemma2}, which directly
  translates into the conditions in (i).

  For the converse, assume (i).  Again, consider the factorization
  $T_Q = T_P^* T_P$ with $T_P$ a lower triangular $3\times 3$ matrix
  with respect to the decomposition in \eqref{eq:11}.  By the
  equivalence of \eqref{eq:4} and \eqref{eq:5} in Lemma \ref{lemma2},
  we have $Y= P_K^*P_K$, where $P_K = \mathrm{row}(P_k)_{k\in K}$.
  Set $P(z) = P_KZ_K = \sum_{k\in K}P_k z^k$.  Then $Q(z) = P(z)^*
  P(z), z \in {\mathbb T}^d$. \end{proof}

The notion of ``outerness'' of the factor $P$ is given above in
equations \eqref{inclusion1} and \eqref{inclusion2}.  These conditions
reduce in the one-variable case to condition (iii) in Theorem
\ref{eqs}.  Clearly, there are many other, perhaps more natural, ways
of generalizing the notion of outerness to the multivariable case
(see, for example, \cite{MR93e:32004}).  For instance, the condition
$\clran T_P ={H^2_{\sp M}({\mathbb D}^d)}$ or the condition that
\begin{equation*}
  P(z)^* P(z) \ge L(z)^* L(z) , z \in {\mathbb D}^d ,
\end{equation*}
for all $L(z)$ with $P(z)^* P(z) = L(z)^* L(z), z \in {\mathbb T}^d$,
are both options.  How all these different notions relate to one
another remains to be investigated.  We leave this for a future
publication.

Recall from \cite{GW} the following result regarding stable
factorization (factorizations in terms of polynomials void of zeros in
$\overline{\mathbb D}^2$) of a strictly positive scalar valued
trigonometric polynomial.

\begin{theorem} \label{gw} \cite{GW}  Let $K = \{ 0, \ldots , n_1\}
  \times \{ 0, \ldots , n_2 \}$ and let $Q_k: {\sp H} \to {\sp H}$, $k
  \in K - K$, be scalar valued so that $Q(z) := \sum_{k \in K-K} Q_k
  z^k > 0, z \in {\mathbb T}^2$.  Then there exists a scalar valued
  polynomial $P(z) = \sum_{k\in K} P_k z^k$ so that $Q(z) = |P(z)|^2$,
  $z \in {\mathbb T}^2$, and $P(z) \neq 0$, $z \in \overline{\mathbb
    D}^2$, if and only if
  \begin{equation*}
    (\Pi_{K\setminus \{ (n_1 , n_2 ) \}} T_{Q^{-1}} 
    \Pi_{K\setminus \{ (n_1 , n_2 ) \}} )^{-1} 
  \end{equation*}
  has zero entries in locations $(k,l)$ where $k \in \{ 1, \ldots ,
  n_1 \} \times \{ 0 \} $ and $l \in \{ 0 \} \times \{ 1, \ldots , n_2
  \}$.
\end{theorem}

The conditions in Theorem \ref{multi} and \ref{gw} are quite
different.  The following theorem, which gives necessary conditions on
the Schur complement in the form of the existence of a decomposition,
relates better to Theorem \ref{gw} as the condition on the Schur
complement implies the necessity of some entries in the Schur
complement being zero.

\begin{theorem} \label{2var}
  Let $K = \{ 0, \ldots , n_1\} \times \{ 0, \ldots , n_2 \}$ and let
  $Q_k: {\sp H} \to {\sp H}$, $k \in K - K$, be Hilbert space
  operators so that $Q(z) := \sum_{k \in K-K} Q_k z^k \ge 0, z \in
  {\mathbb T}^2$.  Put
  \begin{equation*}
    \begin{split}
      S_1 &= S(T_Q; \{ 0, \ldots , n_1-1\} \times \{ 0, \ldots , n_2 \} ), \\
      S_2 &= S(T_Q; \{ 0, \ldots , n_1\} \times \{ 0, \ldots , n_2 -1\} ), \\
      S_0 &= S(T_Q; \{ 0, \ldots , n_1-1\} \times \{ 0, \ldots , n_2 -1\} ). 
    \end{split}
  \end{equation*}
  Suppose that $Q(z) = P(z)^* P(z)$, $z\in {\mathbb T}^2$, where $P(z)
  =\sum_{k \in K} P_k z^k$, $P_k : {\sp H} \to {\sp H}$, satisfies
  \begin{equation} \label{outer2}
    {\rm ran} \Pi_{{\mathbb N}_0^2 \setminus \widetilde K} T_P \Pi_{\widetilde K}
    \subset \overline{\rm ran} \Pi_{{\mathbb N}_0^2 \setminus \widetilde K} 
    T_P \Pi_{{\mathbb N}_0^2 \setminus \widetilde K}, 
  \end{equation}
  for $\widetilde K = \{ 0, \ldots , n_1-1\} \times \{ 0, \ldots , n_2
  -1\} , \{ 0, \ldots , n_1\} \times \{ 0, \ldots , n_2 -1\} , \{ 0,
  \ldots , n_1\} \times \{ 0, \ldots , n_2 -1\} , K \setminus \{ (n_1,
  n_2 ) \}$ and $K$.  Then
  \begin{equation} \label{schureq1}
    S(T_Q; K \setminus \{ (n_1, n_2 ) \} ) = S_1 + S_2 - S_0
  \end{equation}
  and
  \begin{equation}\label{schureq2}
    S(T_Q ; K ) = T_Q - T_1 (T_Q - S_1) T_1^* - T_2 (T_Q - S_2 ) T_2^* +
    T_1 T_2 (T_Q - S_0) T_2^* T_1^* ,
  \end{equation}
  where $T_i$ is the Toeplitz operator corresponding to the
  multiplication operator $M_i :H^2_{\sp H} ({\mathbb D}^2) \to
  H^2_{\sp H} ({\mathbb D}^2)$ with symbol $m_i$, where $m_i(z) =
  z_i$, $i=1,2$.
  
  Conversely, suppose that \eqref{schureq1} and \eqref{schureq2} hold,
  then there exists an operator valued polynomial $P(z) =\sum_{k \in
    K} P_k z^k: {\sp H} \to {\sp H}$ so that $Q(z) = P(z)^* P(z)$,
  $z\in {\mathbb T}^2$ and \eqref{inclusion1} and \eqref{inclusion2}
  hold.
\end{theorem}

\begin{proof} Since $Q(z) = P(z)^* P(z)$, $z \in {\mathbb T}^2$, we
  have that $T_Q = T_P^* T_P$.  Let $\widetilde K = \{ 0,\ldots , p_1
  \} \times \{ 0 , \ldots , p_2 \}$ with $p_i \in \{ n_i, n_i -1 \}$,
  $i=1,2$, or $\widetilde K = K \setminus \{ (n_1,n_2) \}$, and view
  the equation $T_Q = T_P^* T_P$ with respect to the decomposition
  \begin{equation*}
    {\rm ran} \Pi_{\widetilde K} \oplus {\rm ran} 
    \Pi_{{\mathbb N}_0^2 \setminus \widetilde K} . 
  \end{equation*}
  
  Since \eqref{outer2} holds true we have by Lemma \ref{lemma1} that
  \begin{equation*}
    S(T_Q ; \widetilde K ) = \Pi_{\widetilde K} T_P^* \Pi_{\widetilde K} 
    T_P \Pi_{\widetilde K} . 
  \end{equation*}
  This now yields expressions for all operators in \eqref{schureq1}
  and \eqref{schureq2} in terms of $P$.  It is now straightforward to
  check that \eqref{schureq1} and \eqref{schureq2} hold.  For
  illustration purposes let us write out the equalities in the
  operators in case that $n_1=n_2=1$: here we have that
  $S_0=L_0^*L_0$, $S_1=L_1^*L_1$, $S_2=L_2^*L_2$, $S(T_Q; K\setminus
  \{ (1,1) \} ) = L_3^*L_3$, $S(T_Q; K) = L_4^*L_4$, where
  \begin{equation*}
    L_0 =\begin{pmatrix} P_{00} & 0 & 0 & 0 \cr 0 & 0 & 0 & 0 \cr 
      0 & 0 & 0 & 0 \cr 0 & 0 & 0 & 0
    \end{pmatrix} , \qquad 
    L_1 =\begin{pmatrix} P_{00} & 0 & 0 & 0 \cr 
      0 & 0 & 0 & 0 \cr P_{10} & 0 & P_{00} & 0 \cr 0 & 0 & 0 & 0
    \end{pmatrix}, \qquad
    L_2 =\begin{pmatrix} P_{00} & 0 & 0 & 0 \cr 
      P_{01} & P_{00} & 0 & 0 \cr 0 & 0 & 0 & 0 \cr 0 & 0 & 0 & 0
    \end{pmatrix},
  \end{equation*}
  \begin{equation*}
    L_3 =\begin{pmatrix} P_{00} & 0 & 0 & 0 \cr
      P_{01} & P_{00} & 0 & 0 \cr P_{10} & 0 & P_{00} & 0 \cr 0 & 0 & 0 & 0
    \end{pmatrix} , \qquad
    L_4 =\begin{pmatrix} P_0 & 0 & 0 & 0 \cr P_{01} & P_{00} & 0 & 0 \cr 
      P_{10} & 0 & P_{00} & 0 \cr P_{11} & P_{10} & P_{01} & P_{00}
    \end{pmatrix}, 
  \end{equation*}
  and
  \begin{equation}\label{Y0} 
    Y_0 := T_1 T_Q T_1^* - T_2 T_Q T_2^* + T_1 T_2 T_Q T_2^* T_1^*  =
    \begin{pmatrix} Q_{00} & Q_{01}^* & Q_{10}^* & Q_{11}^* \cr 
      Q_{01} & 0 & Q_{1,-1}^* & 0 \cr Q_{10} & Q_{1,-1} & 0 & 0 \cr
      Q_{11} & 0 & 0 & 0
    \end{pmatrix} ,
  \end{equation}
  where we restricted the operators to rows and columns indexed by $\{
  (0,0), (0,1), (1,0), (1,1) \}$, as they contain all the nonzero
  entries.  The operators $T_1$ and $T_2$ restricted to this part
  correspond to
  \begin{equation*}
    \begin{pmatrix} 0 & 0 & 0 & 0 \cr 0 & 0 & 0 & 0 \cr 
      I & 0 & 0 & 0 \cr 0 & I & 0 & 0
    \end{pmatrix}
    \qquad {\rm and } \qquad 
    \begin{pmatrix} 0 & 0 & 0 & 0 \cr I & 0 & 0 & 0 \cr 
      0 & 0 & 0 & 0 \cr 0 & 0 & I & 0
    \end{pmatrix},
  \end{equation*}
  respectively.  Formulas \eqref{schureq1} and \eqref{schureq2} follow
  now directly.  The computations for the case $n_1n_2>1$ are similar.
  
  For the converse we apply Theorem \ref{multi}.  Using
  \eqref{schureq1} and \eqref{schureq2} we find that $Y$ in Theorem
  \ref{multi} equals
  \begin{equation*}
    Y=Y_0 - (S_1 - T_1S_1T_1^*) - (S_2-T_2S_2T_2^*) + (S_0 - T_1T_2 S_0 T_2^*T_1^*) 
  \end{equation*}
  yielding that 
  \begin{equation*}
    Z_K^* Y Z_K = Q(z) - (1-|z_1|^2) (Z_K^* S_1 Z_K) - (1-|z_2|^*) (Z_K^* S_2 Z_K)
    + (1-|z_1z_2|^2) (Z_K^* S_0 Z_K). 
  \end{equation*}
  Thus for $(z_1,z_2) \in {\mathbb T}^2$ we obtain equality
  \eqref{condY}.  The conclusion now follows from Theorem \ref{multi}.
\end{proof}

Notice that Theorem \ref{2var} is not an if and only if statement due
to the different ``outerness'' requirements on $P$: in one direction
the outerness requirement is \eqref{outer2} while in the other
direction it is \eqref{inclusion1} and \eqref{inclusion2}.  We suspect
that these two outerness requirements are different, though we have
not constructed an example showing this.

Note too that \eqref{schureq1} implies that $S(T_Q; K \setminus \{
(n_1 , n_2 )\} )$ has zeros in locations $(k,l)$ where $k \in \{ 1,
\ldots , n_1 \} \times \{ 0 \} $ and $l \in \{ 0 \} \times \{ 1,
\ldots , n_2 \}$.

\bibliographystyle{plain} 

\end{document}